\documentclass[11pt]{article}
\usepackage{mathptmx}

\usepackage{epsfig,amsfonts,latexsym,amssymb,shadow,amsmath,wrapfig}
\usepackage{subfigure}
\usepackage{epstopdf}
\usepackage{verbatim}
\usepackage{mathtools}

\usepackage{graphicx,color,multirow,colortbl,hhline,array,url,fancyhdr,setspace,enumerate}
\usepackage[normalem]{ulem}
\usepackage[margin=1.2in]{geometry}
\usepackage{enumitem}

\usepackage{amsmath,bm}
\usepackage{mathtools} 
\usepackage{float}
\usepackage[ruled,vlined,linesnumbered]{algorithm2e}
\usepackage{booktabs}
\usepackage{tabularx}

\DeclareMathOperator*{\maximize}{maximize}
\usepackage{threeparttable}

\usepackage{bm}

\usepackage{xcolor}

\usepackage[nosort]{cite}

\usepackage[hang,flushmargin,symbol]{footmisc}

\begin{document}
	
	\title{The economics of utility-scale portable energy storage systems in a high-renewable grid}
	
	\begingroup

	\author{Guannan He\footnotemark[2]~~,
		Jeremy Michalek\footnotemark[6]~~\footnotemark[8]~~,
		Soummya Kar\footnotemark[5]~~, 
		Qixin Chen\footnotemark[10]~~,\\ 
		Da Zhang\footnotemark[1]~~\footnotemark[4]~~\footnotemark[3]~~,
		Jay F. Whitacre\footnotemark[1]~~\footnotemark[6]~~\footnotemark[9]
	}

	\footnotetext[1]{Corresponding authors: zhangda@tsinghua.edu.cn, whitacre@andrew.cmu.edu}
	\footnotetext[2]{MIT Energy Initiative, Massachusetts Institute of Technology, Cambridge, MA, USA.}
	\footnotetext[4]{Institute of Energy, Environment and Economy, Tsinghua University, Beijing, China.}
	\footnotetext[3]{MIT Joint Program on the Science and Policy of Global Change, Massachusetts Institute of Technology, Cambridge, MA, USA.}
	\footnotetext[5]{Department of Electrical and Computer Engineering, Carnegie Mellon University, Pittsburgh, PA, USA.}
	\footnotetext[10]{Department of Electrical Engineering, Tsinghua University, Beijing, China.}
	\footnotetext[6]{Department of Engineering and Public Policy, Carnegie Mellon University, Pittsburgh, PA, USA.}
	\footnotetext[8]{Department of Mechanical Engineering, Carnegie Mellon University, Pittsburgh, PA, USA.}
	\footnotetext[9]{Department of Materials Science and Engineering, Carnegie Mellon University, Pittsburgh, PA, USA.}

	\endgroup
	
	\date{}
	
	\maketitle

	\noindent \textbf{Summary}: Battery storage is expected to play a crucial role in the low-carbon transformation of energy systems. The deployment of battery storage in the power gird, however, is currently severely limited by its low economic viability, which results from not only high capital costs but also the lack of flexible and efficient utilization schemes and business models. Making utility-scale battery storage portable through trucking unlocks its capability to provide various on-demand services. We introduce the potential applications of utility-scale portable energy storage and investigate its economics in California using a spatiotemporal decision model that determines the optimal operation and transportation schedules of portable storage. We show that mobilizing energy storage can increase its life-cycle revenues by 70\% in some areas and improve renewable energy integration by relieving local transmission congestion. The life-cycle revenue of spatiotemporal arbitrage can fully compensate for the costs of portable energy storage system in several regions in California, including San Diego and the San Francisco Bay Area. 
	
	%
	
	\section{Introduction}
	\label{sec:studythree:intro}
	Energy storage will be essential in future low-carbon energy systems to provide flexibility for accommodating high penetrations of intermittent renewable energy \cite{chu12,ALBERTUS202021,braff16,stephan16}. Currently, the scale of existing utility-scale battery energy storage capacity is still relatively low compared to installed wind and solar capacities as the return of energy storage investment is inadequate due to the high upfront costs and the lack of flexible and efficient schemes for storage utilization \cite{HAMELINK2019120,LOMBARDI2017485}. While demands for flexibility (such as time shift \cite{davies2019combined}, congestion relief \cite{del14}, and ramping \cite{CUI201827}) supplied by energy storage will become increasingly pervasive \cite{fares2017impacts,cochran2014flexibility}, they are intermittent and distributed---varying across both time and location---and thus usually result in a low utilization rate if the energy storage system is deployed at a fixed location. For example, in a time-shift application, the energy storage system will operate only when electricity prices reach extremes as a result of very high or low renewable generation and/or electricity demand and stay idle most of the time \cite{he18}. Similar low-utilization patterns are observed for grid congestion relief applications \cite{elliott18,lo2012}, and flexible ramping capacities are required only when there is a significant fluctuation in renewable energy or demand (e.g., at the sunrise and sunset) \cite{wang17enhancing}. In addition to economic disadvantages of intermittent revenue streams, low utilization rates also shorten the revenue-generating lifetime of battery storage system due to calendar degradation \cite{wang14,xu18,ecker14} and undermine its economic viability.
	
	\paragraph{}
	Better use of storage systems is possible and potentially lucrative in some locations if the devices are portable, thus allowing them to be transported and shared to meet spatiotemporally varying demands \cite{elliott18}. Existing studies have explored the benefits of coordinated electric vehicle (EV) charging \cite{Weis15,wolinetz18}, vehicle-to-grid (V2G) applications for EVs \cite{shao17,alizadeh17} and railway systems \cite{sun15,sun17-iHtRb} as well as EVs supplying capacities for emergency scenarios in power distribution systems \cite{kim18}. Routing problems for EVs with V2G option have also been studied, though with limited temporal resolution \cite{cabrera19} or decision flexibility \cite{abdulaal17,tang19}. To the best of our knowledge, there is no existing work that systematically investigates the potential applications and related economics of utility-scale portable energy storage using a comprehensive spatiotemporal decision model.
	
	\paragraph{}
	In this work, we first introduce the concept of utility-scale portable energy storage systems (PESS) and discuss the economics of a practical design that consists of an electric truck, energy storage, and necessary energy conversion systems. In this business model, the truck is loaded with energy storage and travels to provide on-demand services within a certain area. We develop a spatiotemporal decision model that determines the optimal operation and transportation schedules of the PESS to maximize its profit. The model is applied to study the economics of this PESS design in California, in which the PESS helps integrate renewable energy and relieve grid congestion at the same time. We find that compared to the stationary energy storage system (SESS), the life-cycle revenue of PESS can be 70\% higher in some areas. In fact, the spatiotemporal arbitrage could generate revenue high enough to recover the upfront cost of the storage system and becomes one of the most profitable grid applications for utility-scale energy storage in California.
	
	\section{Portable energy storage system}
	\label{sec:studythree:portable_storage}
\begin{figure}[!ht]
		\centering
		\includegraphics[width=\textwidth]{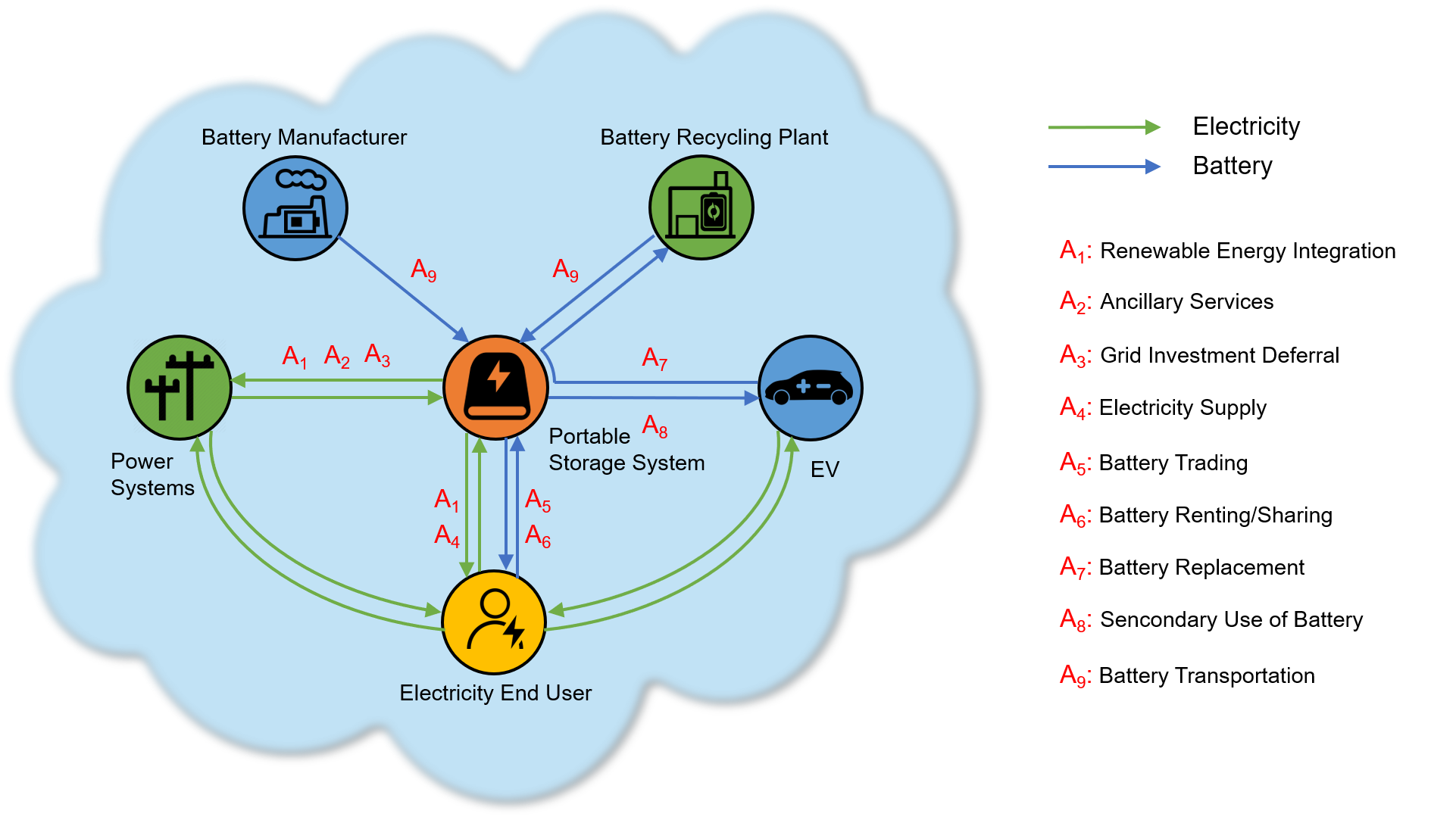}
		\caption{Schematic of energy and battery flows and potential applications of a PESS in energy and transportation systems. The green arrows indicate that the PESS exchanges energy with the entities to which it provides services, while the blue arrows represent battery flows/exchanges between the PESS and other entities.}
		\label{fig-1}
	\end{figure}

	A typical PESS integrates utility-scale energy storage (e.g., battery packs), energy conversion systems, and vehicles (e.g., trucks, trains, or even ships). The PESS has a variety of potential applications in energy and transportation systems and can switch among different applications across space and time serving different entities, like a cloud of on-demand resource, as shown in Figure \ref{fig-1}. PESSs can provide the same services as SESSs, such as renewable energy integration, various ancillary services, grid congestion relief to defer investments, and so on. But the portability of PESS also enhances its capability to tap into multiple value streams that have spatiotemporal variability, which in turn improves its asset utilization and potentially its value proposition over the SESS. When renewable energy integration is limited by grid transmission capacity, a PESS taking advantage of spatiotemporal arbitrage opportunities by traveling between grid nodes with congestion (where constructing new transmission lines is cost-inefficient and time-consuming) to charge at low-price nodes with overabundant renewable energy and discharge at high-price nodes can integrate more renewable energy and thus generate a higher revenue than a SESS. 
	
	\paragraph{}
	Besides spatiotemporal arbitrage, PESSs can also provide an effective way to cope with seasonal short-term power shortages, which are likely to become more frequent in the future due to extreme weather events caused by climate change or high penetration of intermittent renewable supply in deep decarbonization scenarios \cite{ZIEGLER2019}. The huge capacities of energy storage required to avoid power outage in some short periods can potentially be borrowed and transported from other regions using PESSs. PESSs can also serve as physical platforms for battery trading, sharing, and reuse, complemented by on-demand financial contracts \cite{thomas2019general,wang2018transactive}. In addition, PESSs can potentially support recycling and reuse of batteries from EVs, saving battery transportation costs in the cycles. EV users will have the option to replace used batteries with new ones at lower prices based on the states of health of the original batteries. Electricity consumers could also potentially rent or replace batteries from PESSs and reduce demand charges.

	\paragraph{}
	The PESS studied in this paper involves loading lithium-ion batteries and inverters onto containers and trailers that can be hauled by electric trucks, supplemented with battery, thermal, and energy management systems for safety and control purposes. We choose Tesla Powerpack and Tesla Semi as the battery and vehicle with basic technical and cost specifications estimated in Table \ref{tab-1}. The price of a Tesla Semi with a 500-mile range is about \$180,000, which should include an approximately 1-MWh battery according to its energy consumption rate and range \cite{Sripad19}. Based on the payload of Tesla Semi and the densities of Tesla Powerpack and inverter, one truck could accommodate approximately 2.7 MWh batteries with inverters. The total capital cost of a PESS is approximately \$735,000, assuming a unit cost of \$200/kWh for battery packs.
	
	\begin{table}[!htb]
		\centering	
		\begin{threeparttable}
			\renewcommand{\arraystretch}{1.5}
			\caption{Estimates of the technical and cost specifications for a PESS design consisting of Tesla Semi and Powerpack.}
			\label{tab-1}	
			\begin{tabular}{l c}
				\hline\hline
				Capital cost of Tesla Semi with 500-mile range (US \$) \tnote{1} & 
				180,000 \\
				\midrule
				Energy consumption rate of Tesla Semi (kWh/mile) \tnote{1} & 
				$<$2  \\
				\midrule
				Tesla Semi freight payload (tonne) \cite{Sripad19}& 
				19\\
				\midrule
				Energy density of Tesla Powerpack with inverter (kWh/kg) \tnote{2}  & 
				0.11  \\
				\midrule
				Power density of Tesla Powerpack inverter (kVA/kg) \tnote{2}  & 
				0.63  \\
				\midrule
				Total battery energy capacity per truck (MWh, 1-hour duration) & 
				2.7 \\
				\midrule
				Capital cost of Tesla Powerpack (US \$, given \$200/kWh \cite{inverter18})& 
				340,000 \\
				\midrule
				Capital cost of inverter (US \$, given \$70/kW \cite{inverter18})& 
				190,000 \\
				\midrule
				Capital cost of trailer (US \$) \cite{moultak2017transitioning}& 
				25,000 \\	
				\midrule
				Total capital cost of a PESS (US \$)& 
				735,000 \\	
				\hline\hline
			\end{tabular}
			\begin{tablenotes}
				\item[1] Source: https://www.tesla.com/semi
				\item[2] Source: https://www.tesla.com/powerpack 
			\end{tablenotes}
		\end{threeparttable}
	\end{table}
	
	\section{Spatiotemporal arbitrage revenue of PESS in California}
	\label{sec:studythree:case_study}
	Here we evaluate the spatiotemporal arbitrage revenues of a PESS in California, where intensive and extensive local grid transmission congestion has been observed recently (see Figure S1). We applied a spatiotemporal arbitrage optimization model (see Methods) to a PESS operating over 1,131 case areas in California. Each case area has a 10-mile radius and takes one of the grid nodes defined by the California Independent System Operator (CAISO) as its center. CAISO also publishes the wholesale locational marginal prices (LMPs) for each node based on the cost of generating and delivering electricity to the node. The number of grid nodes in a case area ranges from 1 to 44, depending on the node density (usually correlated to population density). The PESS is assumed to be a price-taker, which means that the actions of the PESS have a negligible impact on the prices in the case areas. We also assume that the PESS can provide reserve as a secondary application. Based on the LMPs and non-spinning reserve prices in CAISO day-ahead markets, we optimize the operation and transportation strategies of a PESS for each day in its lifetime using the spatiotemporal arbitrage model and calculate its life-cycle revenue as the sum of discounted daily revenues.
	
\begin{figure}[!htb]
		\centering
		\includegraphics[width=\textwidth]{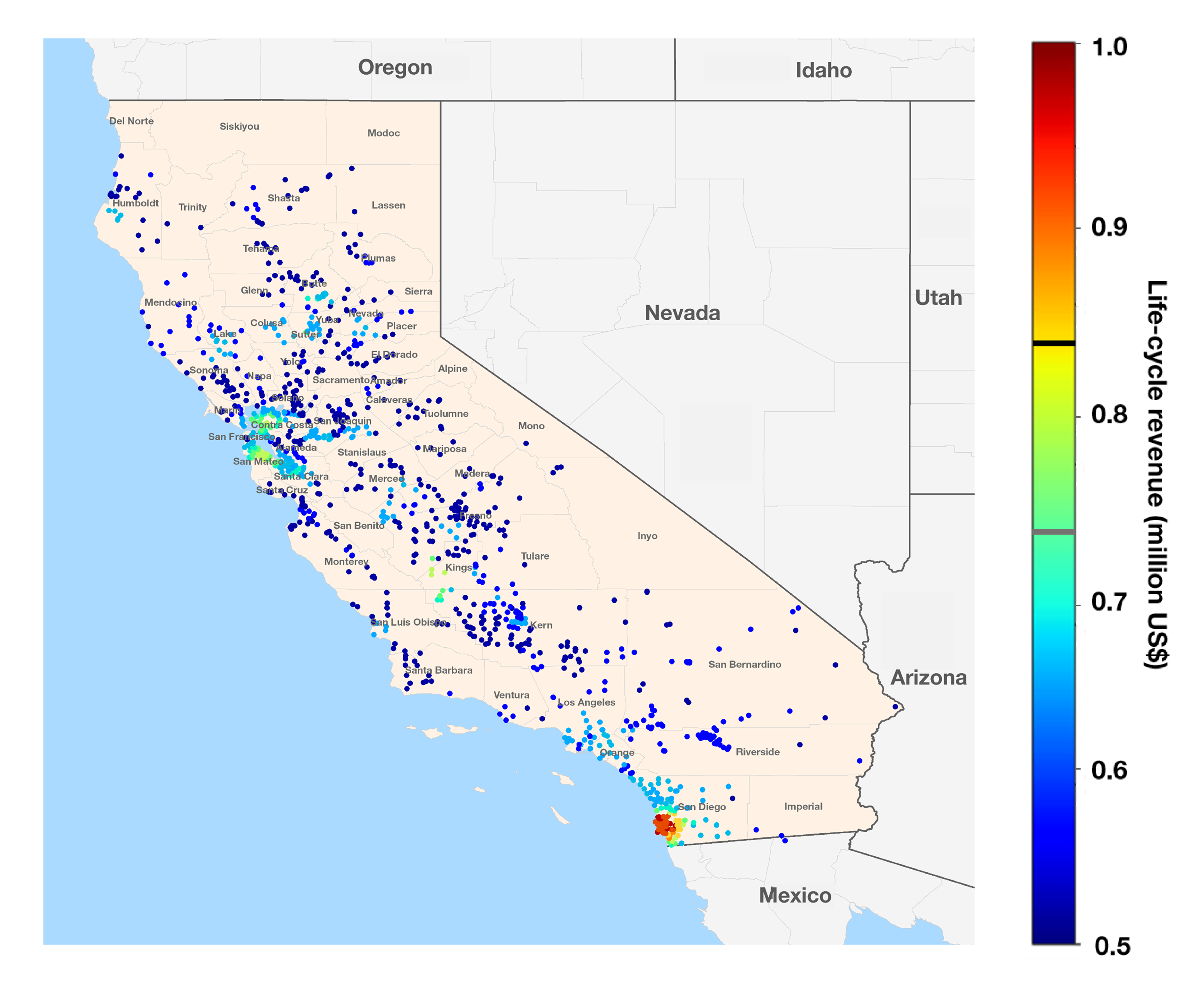}
		\caption{Revenue potential of the 2.7-MWh PESS setup when executing optimal spatiotemporal energy arbitrage in California. Each dot represents a case in which the PESS operates over an area with the dot as its center and a 10-mile radius. The dot color indicates the magnitude of life-cycle revenue for each case. The up-front capital costs plus the fixed operation and maintenance (O\&M) costs (2\% of the capital cost annually \cite{beuse20}) are indicated in the color bar by the black line given a \$200/kWh unit cost for battery packs and the grey line given a \$150/kWh unit cost for battery packs, respectively.}
		\label{fig-2}
	\end{figure}

	\paragraph{}
	Figure \ref{fig-2} presents the resulting life-cycle revenues of PESSs for each case area in California. Each dot represents a case, and the dot color indicates the amount of life-cycle revenue. High-revenue opportunities exist around metropolitan areas, including San Diego, Los Angeles, and the San Francisco Bay Area. The high population densities in these areas correspond to relatively high peak electricity demands and also high penetrations of roof solar generation, which increase the volatility of electricity market prices and create arbitrage opportunities for the PESS. Moreover, the frequency of grid congestion around metropolitan areas is also higher because transmission capacity expansion is more costly and time-consuming, offering a unique advantage to the PESS for grid congestion relief and deferring grid investments. It is important to note that in some areas with lower population densities, such as Kings County and Sutter County, there are also favorable spatiotemporal arbitrage opportunities due to overabundant solar energy from large solar farms.
	
	\begin{figure}[!htb]
		\centering
		\includegraphics[width=\textwidth]{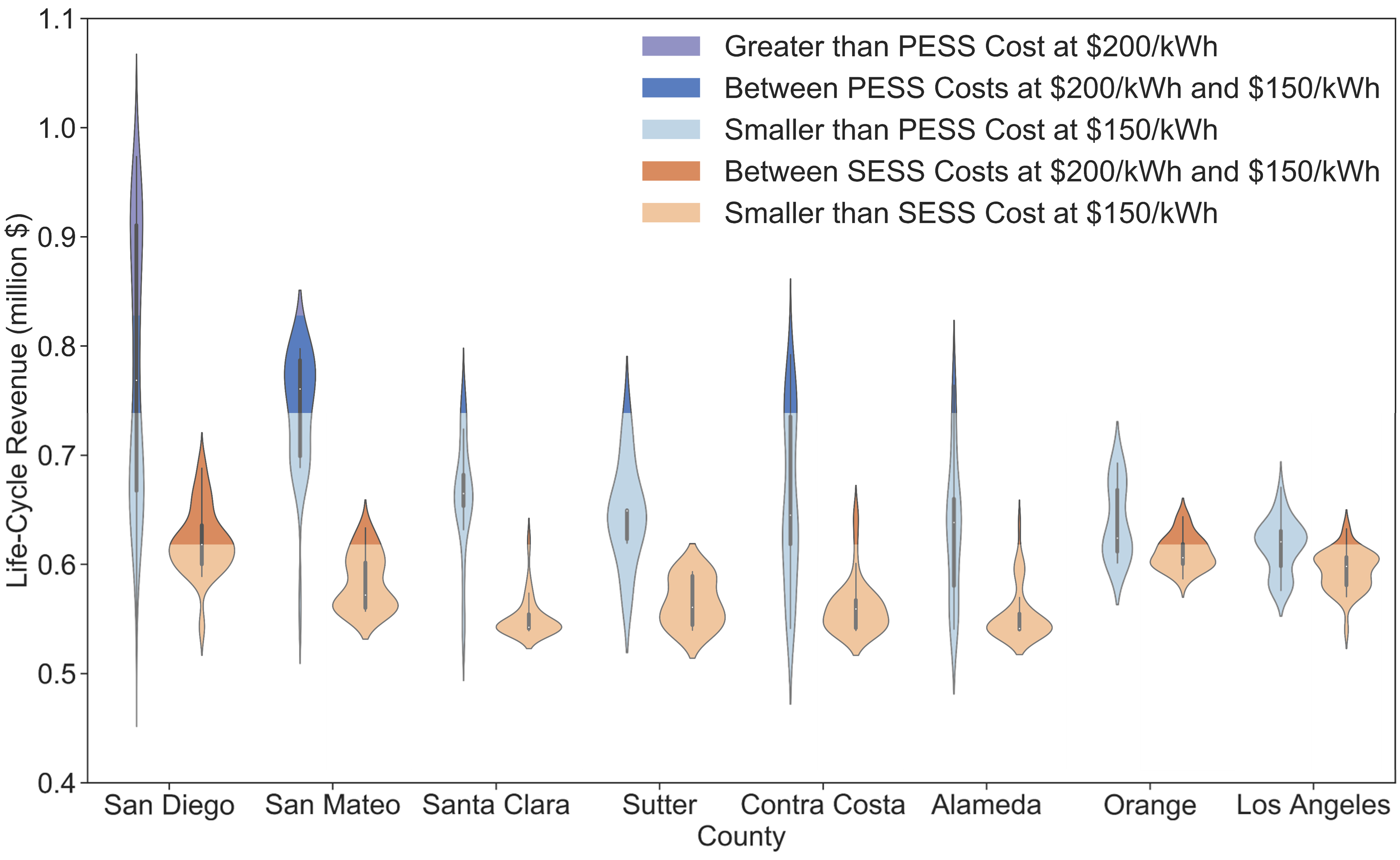}
		\caption{Life-cycle revenues of PESS and SESS in California by county. The top eight counties with the highest median PESS revenue and at least ten case areas are displayed and sorted in descending order of median PESS revenue. Each violin plot shows the median (the white dot), the first and third quartiles (the black bar), the upper and lower adjacent values (the black line), and the distribution of life-cycle revenues for all the case areas in each county. The shading colors of each violin plot represent whether the life-cycle revenue is greater or smaller than the total fixed costs (up-front capital costs plus the fixed O\&M costs) of PESS/SESS given a \$200/kWh or \$150/kWh unit cost for battery packs. The total fixed cost of a SESS is approximately \$0.13 million lower than that of a PESS to exclude the extra capital costs from the powertrain, trailer, etc. and the extra O\&M costs from insurance, licensing, etc. \cite{Sripad19,moultak2017transitioning,beuse20}}
		\label{fig-3}
	\end{figure}
	\paragraph{}
	The life-cycle revenues of PESS and SESS are compared in Figure \ref{fig-3} across all case areas in each of the eight counties in California with the highest median revenues for a PESS. We found significant revenue improvements from SESS to an identically-sized PESS in most of the counties, with the highest increase of approximately 70\% for the case area in San Diego with the highest life-cycle revenue. The median revenue across case areas for a PESS is 25\% (\$0.15 million) greater than the median revenue of a SESS in San Diego County, 27\% (\$0.19 million) in San Mateo County, and 21\% (\$0.12 million) in Santa Clara County. Moreover, a PESS can generate enough revenues to fully compensate for its total costs in some case areas, while a SESS cannot in most case areas. For a PESS, the median revenues of San Diego and San Mateo counties are higher than the total cost given the near-term cost estimate (\$150/kWh for battery pack). Sometimes, even travelling and arbitraging between only two nodes could produce enough revenues to cover the capital cost in some areas. For example, in the case of two nodes around Kettleman City in Kings County, the life-cycle revenue of PESS is \$0.76 million, surpassing the near-term cost. From the distributions of the difference in estimated profitability between PESS and SESS across case areas in the eight counties (see Figure S8), we also found that the PESS is more profitable than the SESS in cases spanning over 36\% of all 33 studied counties with at least ten case areas in California. The above results indicate that converting SESS to PESS can potentially make turnrounds in profitability.

	\section{Operational patterns of a PESS in spatiotemporal arbitrage}
\begin{figure}[!ht]
		\centering
		\includegraphics[width=\textwidth]{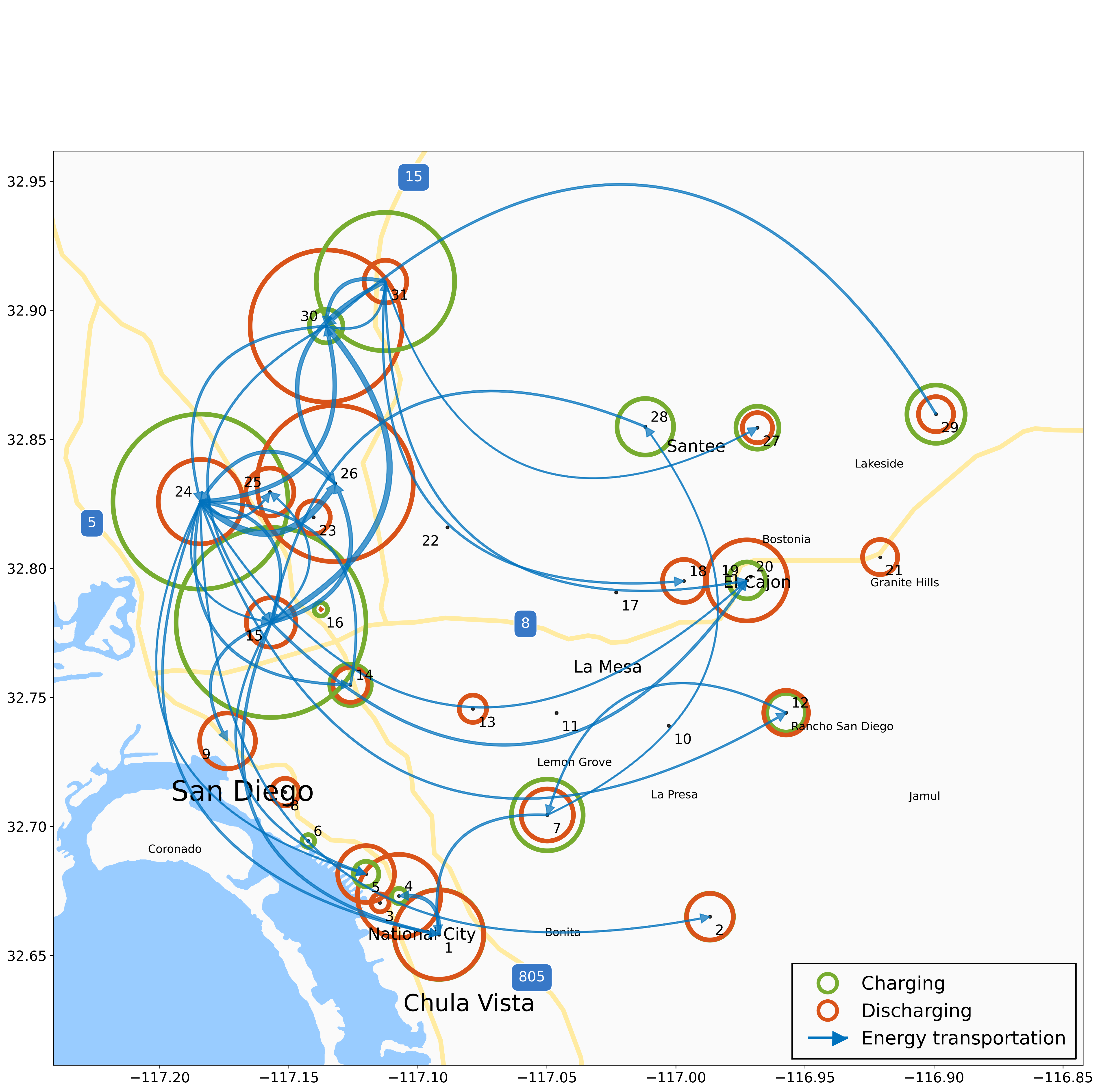}
		\caption{Optimal charging, discharging, and traveling strategies of the PESS for spatiotemporal arbitrage in a case area in San Diego County, California for one year. The case area has a 10-mile radius, with node 17 (CAISO node ID: MURRAY\_6\_N005) as its center. The sizes of the red and green circles represent the amounts of charged and discharged energy, respectively, for the node at the center of the circles. The arrows represent the directions of energy transmission by the PESS from one node to the other. The line width of the arrows indicates the amount of transmitted energy.}
		\label{fig-4}
	\end{figure}

	To further reveal how the profitability of a SESS can be enhanced by mobilization, we show the optimal operational strategies of a PESS in the case area with the highest life-cycle revenue (\$0.97 million) in California located in San Diego County. In Figure \ref{fig-4}, for each node, the circles around it represent the amounts of energy that the PESS charges from (green circle) and discharges to (red circle) the node in a year; and the directions and the line widths of the arrows represent the directions and the amount of energy transmission by the PESS between nodes. We can observe that the PESS frequently operates at some nodes, including node 15, 24, 26, 30, and 31 in Figure \ref{fig-4}. Some nodes both import and export energy, such as node 7, 15, 24, 31, and etc. The most frequently travelled route is between node 24 and 26 (route 24-26), which are also the nodes with the most frequent charging and discharging, respectively. Although route 24-26 is comparatively short, there are also many less frequently travelled long routes, with an aggregate energy transmission comparable to that of route 24-26. In summary, profit opportunities exist in the area and are widely distributed, thus enabling energy storage to be shared across different nodes is critical to fully exploit these opportunities. 
	\begin{figure}[!ht]
		\centering
		\includegraphics[width=\textwidth]{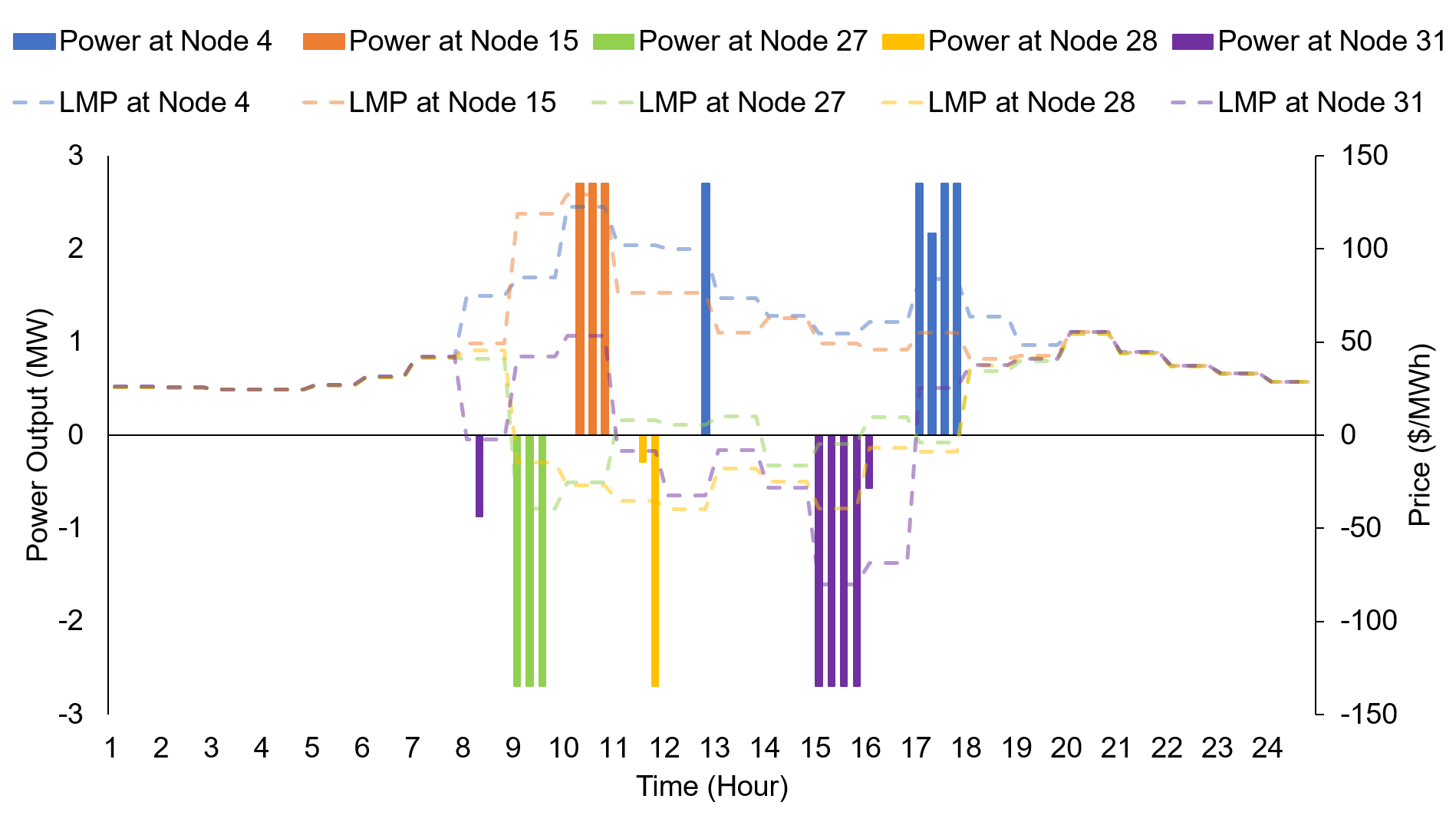}
		\caption{Optimal operational schedules for a PESS doing spatiotemporal arbitrage in a case area in San Diego County and the locational marginal price (LMP) profiles of the relevant nodes in a sample day (Mar 23, 2018). The bars indicate the charging (below the x-axis) and discharging (above the x-axis) power, and the dashed lines represent the LMP profiles.}
		\label{fig-5}
	\end{figure}
	
	\paragraph{}
	We choose the day with the most frequent travel of the year as a sample day to show the optimal daily operational schedule of the PESS in Figure \ref{fig-5}. The PESS visits five nodes in the sample day, and the dashed lines show the locational marginal prices (LMPs) of the five nodes. From 8 am to 8 pm, there are remarkable price differences among the five nodes, indicating grid transmission congestion. Node 4 and 15 are high-price nodes and node 27, 28, and 31 are low-price nodes. The negative prices at low-price nodes in some periods imply excess solar generations. To exploit the price differences, the PESS travels among the five nodes to charge and buy energy at the node with the lowest price and discharge and sell energy at the node with the highest price. In Figure \ref{fig-5}, the PESS makes a trip from one node to the other at the time represented by gaps between bars with different colors. Five trips are made in the sample day. The traveling capability provides the PESS much more profit opportunities compared to the SESS, because the PESS can profit from both price differences between different nodes and between different hours within one node, while the stationary storage can only profit from the price difference within one node. As seen from Figure \ref{fig-5}, the PESS conducts three profitable charge-discharge cycles over the day, while for the SESS it is usually only profitable to run one cycle during the peak and valley hours (see  Figure S9). If we use the amount of energy charged to storage at negative LMPs to approximate the amount of integrated renewable energy by storage, the PESS accommodates more than four times the renewable energy as the SESS does in a year (238 MWh versus 58 MWh). The transmission congestion in the area can also be relieved at the same time.
	
	\clearpage
	\section{Discussions}
	\label{sec:studythree:conclusions}
	We introduce and assess a new business model for energy storage deployment in which battery packs are mobilized to provide various types of on-demand services in energy and transportation systems. The new portable deployment has many potential applications that a stationary deployment cannot be used for. We develop a spatiotemporal optimization model for a PESS and apply the model to simulate the operation and transportation of a PESS design consisting of a Tesla Semi and Powerpack in California to perform arbitrage within a 10-mile radius of each of 1131 case sites. We find that portable deployment has the potential to enhance profitability relative to stationary deployment in 36\% of the studied counties and to exceed costs in San Diego sites as well as several other locations as battery costs drop.
	
	\paragraph{}
	While frequency regulation has been recognized as one of the most profitable grid application for utility-scale energy storage in many regions of the world \cite{eyer2010energy,RAPPAPORT2017609}, the average benefit of usage (the life-cycle revenue divided by the total available energy throughput) for spatiotemporal arbitrage is approximately \$70/MWh-throughput (San Diego County median value), more than twice of that for frequency regulation in California \cite{he18}. 
	The daily revenue for spatiotemporal arbitrage is approximately \$43/MWh-capacity, also higher than the values for any other grid application reported in \cite{davies2019combined}.
	
	\paragraph{}
	In addition to the profitability improvement of energy storage, transmission congestion can be relieved as the PESS transmits energy across nodes, and the transmission investment can be saved or delayed. The PESS has several advantages over transmission capacity expansion. First, it can be shared among multiple congestion areas and thus has higher potential utilization rates than a new transmission line; second, the PESS can be committed much faster when new congestion emerges due to distributed energy resource integration, while building new transmission can take over ten years; third, a PESS can move within a power system or across different systems to adapt to seasonal or longer-term changes in renewable resources and demands as systems evolve and the climate changes \cite{Zeng19,Karnauskas17,DOWLING2020}. It is worth mentioning that the PESS will not replace, but only complement, transmission lines. How to optimally co-plan PESS and transmission lines remains an interesting question to be explored.
	
	\paragraph{}
	Some factors may affect the PESS profitability either positively or negatively. We limit the range of each case area to a circle with a 10-mile radius to keep the spatiotemporal decision model computationally tractable, which results in profit underestimation. Imperfect price forecasting and traffic congestion can both reduce the PESS revenue, which we further discuss in the Supplemental Notes. We did not consider the value proposition changes due to infrastructure changes in power system over the PESS lifetime either, which can be either negative or positive. If the price differences in a region are shrinking, the PESS can either travel to other regions or switch to other applications to alleviate this negative impact of value proposition change--a key advantage of PESSs over SESSs. 
	
	\paragraph{}
	As introduced in Figure \ref{fig-1}, there are many other potential applications of PESSs besides energy arbitrage in future energy and transportation systems. The value propositions of some applications are hard to accurately assess due to lack of markets, such as serving as a platform for battery renting/sharing, promoting battery secondary use, and so on. As markets emerge with the increasing penetrations of EV and residential PV, the spatiotemporal decision model proposed in this paper can also be applied to assess the value of those applications and support the optimal dispatch of PESS. 
	
	\paragraph{}
	Currently, there are some policy barriers that limit energy storage to combine different value streams, for both PESS and SESS \cite{stephan16}. Another potential issue is that whether the storage can be assumed as a price-taker or, in other words, how to compensate the storage if it resolves the congestion and eliminates the price difference, where the arbitrage revenue comes from. Addressing this issue calls for policy or market innovations on mechanisms for storage compensation and grid congestion information release. Some technical issues such as safety and thermal control may also need attention in practice. There may be concerns about battery explosion on the road, although its risk may be lower than trucks full of gasoline, whose energy density is an order of magnitude higher than that of lithium-ion batteries.
	 
	\section{Methods}
	
	\subsection{Spatiotemporal Decision Model}
	A spatiotemporal decision model is developed for a PESS to maximize its profit in a region subject to operation and transportation constraints.
	
	\subsubsection*{Objective function}
	The objective $ f(\mathbf{x}_t) $ of the spatiotemporal decision model when optimizing time period $ t $ (typically one day) is to maximize the total market revenue of the portable storage $ R(\mathbf{x}_t) $ minus the transportation cost $ C^{\mathrm{TRA}}(\mathbf{x}_t) $ and the degradation cost $ C^{\mathrm{DEG}}(\mathbf{x}_t) $, as in equation (\ref{4.1}). 
	
	\begin{equation}
	\begin{split}
	\label{4.1}
	&\maximize\limits_{\mathbf{x}_t} f(\mathbf{x}_t)  =
	R(\mathbf{x}_t) - C^{\mathrm{TRA}}(\mathbf{x}_t) - C^{\mathrm{DEG}}(\mathbf{x}_t) \\
	&\mathrm{where}\ \mathbf{x}_t = [P^{\mathrm{DIS}}_{nh},P^{\mathrm{CHA}}_{nh}, \gamma_{nmh}, \omega_{nh}, \alpha_{nh}, \beta_{nh}, \theta_{nh} \quad \forall n,m\in \Omega, h\in \mathcal{H}_t]
	\end{split}
	\end{equation}
	
	The decision variables $\mathbf{x}_t$ include schedules for discharging $ P^{\mathrm{DIS}}_{nh}\in\mathbb{R_+} $, charging $ P^{\mathrm{CHA}}_{nh}\in\mathbb{R_+} $, travel between locations $ \gamma_{nmh}\in\{0,1\} $, parking at location $\omega_{nh}\in\{0,1\}$, arrival to location $\alpha_{nh}\in\{0,1\}$, departure from location $\beta_{nh}\in\{0,1\}$, and an auxiliary location variable $\theta_{nh}\in\{0,1\}$ for the PESS across nodes and time, where $n,m\in\Omega$ indexes the set of grid access nodes where charging and discharging can occur within the specified radius of the case site and $ h\in\mathcal{H}_t $ indexes the set of time intervals, each of length $\Delta h$, within day $t$ during which charging, discharging, and travel decisions are made. We use $\Delta h = 15$ minute intervals in our case studies. For the travel variables, $ \gamma_{nmh} $ indicates whether or not the PESS is traveling from node $ n $ to node $ m $ during time interval $ h $; $\omega_{nh}$ indicates whether the PESS is parked at node $n$ during time interval $h$; $\alpha_{nh}$ indicates whether the PESS arrives at node $n$ during time interval $h$; $\beta_{nh}$ indicates whether the PESS departs from node $n$ during time interval $h$; and $\theta_{nh}$ is a dummy variable used to ensure consistency between $\alpha_{nh}$ and $\gamma_{nh}$.
	 
    \paragraph{}
    The PESS revenue is expressed in equation (\ref{4.2}), where $ \lambda_{nh} $ is the LMP at node $ n $ and time $ h $. For PESS applications other than energy arbitrage, $ \lambda_{nh} $ can be replaced by any benefit rate received by the PESS. The impact of LMP forecasting error on PESS revenue in spatiotemporal arbitrage is evaluated in Figure S4. 
	
	\begin{equation}
	\begin{split}
	\label{4.2}
	R(\mathbf{x}_t) = \Delta h\sum\limits_{h \in \mathcal{H}_t} \sum\limits_{n \in \Omega}
	\lambda_{nh} \left( P^{\mathrm{DIS}}_{nh} - P^{\mathrm{CHA}}_{nh}
	\right)  
	\end{split}
	\end{equation}
	
	\paragraph{}
	The main transportation cost is labor cost, which is assumed to be proportional to the total travel time during day $ t $, as in equation (\ref{4.3}), where $ c_{\mathrm{TRA}} $ denotes the transportation cost per unit time. We use a \$20/hour labor cost in the case studies. The energy consumption during transportation is less than 2 kWh/mile, which translates to 50 kWh/hour given a 25 mile/hour speed. Considering that the PESS always charges at low prices, e.g., below \$20/MWh, the cost of transportation energy consumption is less than \$1/hour. So $ c^{\mathrm{TRA}} $ is set to \$20/hour.
	
	\begin{equation}
	\begin{split}
	\label{4.3}
	C^{\mathrm{TRA}}(\mathbf{x}_t) = c^{\mathrm{TRA}} \Delta h
	\sum\limits_{h \in \mathcal{H}_t} \sum\limits_{n \in \Omega}
	\sum\limits_{m \in \Omega}
	\gamma_{nmh}   
	\end{split}
	\end{equation}
	
	\paragraph{}
	Equation (\ref{4.4}) presents the degradation cost of portable storage. The marginal degradation cost coefficient $c_{t}^{\mathrm{DEG}}$ reflects the opportunity cost of battery usage and is equal to a constant divided by a discount factor $c^{\mathrm{DEG}}/\delta_t$. We use a typical exponential discounting: $ \delta_{t} = (1+r)^{-\kappa(t)} $, where $ r $ is the discount rate (7\% in this study), and $ \kappa(t) $ is the year number for time $ t $ from the beginning of the battery project. The life-cycle marginal degradation cost $c^{\mathrm{DEG}}$ is set to \$50/MWh-throughput, which is determined using an intertemporal decision framework \cite{he18} to achieve the maximum life-cycle revenue for a PESS in the spatiotemporal arbitrage application. $ q_{t} $ is the calendar degradation of the PESS, 1 MWh-throughput/day in the case studies, which is translated from approximately 1\% capacity loss per year \cite{ecker14,grolleau14,keil16}. The cycling degradation is typically a function of battery charging profile. An approximate cycling degradation function for energy arbitrage applications \cite{he18} is expressed in the term following calendar degradation. 
	\begin{equation}
	\begin{split}
	\label{4.4}
	& C^{\mathrm{DEG}}(\mathbf{x}_t) = c_{t}^{\mathrm{DEG}}
	u_{t}\\
	& \mathrm{where}\ u_{t} = 
	q_{t} + \sum\limits_{h \in \mathcal{H}_t} 
	\left( P^{\mathrm{DIS}}_{nh} + P^{\mathrm{CHA}}_{nh}
	\right) \Delta h
	\end{split}
	\end{equation}
	
	\paragraph{}
    It should be noted that the degradation cost is not a real cost but an opportunity cost and thus should be added back to the objective to calculate the real maximum profit as
    $f(\mathbf{x}^{\ast}_{t}) +  C_{t}^{\mathrm{DEG}}(\mathbf{x}^{\ast}_{t})$
    \cite{he18}, where $\mathbf{x}^{\ast}_{t}$ is the optimizer of equation (\ref{4.1}).
	
	\subsubsection*{Storage operation constraints}
	The energy constraints of storage are formulated in equation (\ref{4.6}). The energy level of storage at time $ h $, $ E_{h} $, is a function of the energy level at time $ h-1 $ and the charging/discharging schedules at time $ h $, where $ \rho $ is the self-discharge rate, and $ \eta $ is the charge/discharge efficiency. We set $ \rho $ to 0 and $ \eta $ to 95\% in our case studies. The energy level of storage cannot exceed its capacity, $ E^{\mathrm{MAX}} $ or drop below zero.

	\begin{equation}
	\begin{split}
	\label{4.6}
	&0 \leq E_{h} \leq E^{\mathrm{MAX}}
	\qquad \forall h \in \mathcal{H}_t \\
	&\mathrm{where}\ E_{h} = (1-\rho)E_{h-1}
	+ \sum\limits_{n \in \Omega}
	\left( 
	P^{\mathrm{CHA}}_{nh}\eta\Delta h
	- P^{\mathrm{DIS}}_{nh}\Delta h/\eta
	\right) 
	\end{split}
	\end{equation}
	
	\paragraph{}
	The power output constraints are expressed as equation (\ref{4.7DIS}) and equation (\ref{4.7CHG}), where $ P^{\mathrm{MAX}} $ is the power capacity of the storage, and $ \omega_{nh}\in\{0,1\} $ is a binary variable that denotes whether the storage is at node $ n $ during time interval $ h $ (1 indicates present and 0 indicates absent), a location indicator. This indicator couples the operation and transportation constraints.
	
	\begin{equation}
	\begin{split}
	\label{4.7DIS}
	0 \leq P^{\mathrm{DIS}}_{nh} \leq \omega_{nh} P^{\mathrm{MAX}}
	\qquad \forall n \in \Omega, h \in \mathcal{H}_t
	\end{split}
	\end{equation}
	\begin{equation}
	\begin{split}
	\label{4.7CHG}
	0 \leq P^{\mathrm{CHA}}_{nh} \leq \omega_{nh} P^{\mathrm{MAX}}
	\qquad \forall n \in \Omega, h \in \mathcal{H}_t
	\end{split}
	\end{equation}
	
	\subsubsection*{Storage transportation constraints}
	The storage can only be present at one node at one time and cannot be parked at a node when it is traveling between nodes: 
	\begin{equation}
	\begin{split}
	\label{4.8}
	\sum\limits_{n \in \Omega} \omega_{nh}
	\leq
	1 -  \sum\limits_{n \in \Omega}\sum\limits_{m \in \Omega}
	\gamma_{nmh}
	\qquad \forall h \in \mathcal{H}_t
	\end{split}
	\end{equation}
	
	\paragraph{}
	The traveling status of storage is modelled in equations (\ref{4.9})-(\ref{4.13}), where $ \alpha_{nh}\in\{0,1\} $ is a binary arrival variable that denotes whether the PESS is traveling to node $ n $ at time $ h $; $ \beta_{nh}\in\{0,1\} $ is a binary departure variable that denotes whether the PESS is traveling from node $ n $ at time $ h $; and $ \theta_{nh}\in\{0,1\} $ is an auxiliary binary variable. Specifically, equation (\ref{4.9}) enforces that the arrival indicators $\alpha_{nh}$ and departure indicators $\beta_{nh}$ are consistent with changes in the location indicators $\omega_{nh}$; equation (\ref{4.10}) ensures that arrival and departure are not simultaneously indicated at the same time and place; equation (\ref{4.11}) enforces that travel from a node is indicated once departure from the node is indicated; equations (\ref{4.12})-(\ref{4.13}) ensure that arrival indicators are equal to 1 in time intervals where travel to the node changes from 1 to 0 and equal to 0 otherwise. 
	
	\begin{equation}
	\begin{split}
	\label{4.9}
	\alpha_{nh} - \beta_{nh}
	= \omega_{nh} - \omega_{n(h-1)}
	\qquad \forall n \in \Omega, h \in \mathcal{H}_t
	\end{split}
	\end{equation}
	\begin{equation}
	\begin{split}
	\label{4.10}
	\sum\limits_{n \in \Omega}\left( \alpha_{nh} + \beta_{nh}
	\right) 
	\leq 1 
	\qquad \forall h \in \mathcal{H}_t
	\end{split}
	\end{equation}
	\begin{equation}
	\begin{split}
	\label{4.11}
	\sum\limits_{m \in \Omega}\gamma_{nmh}
	\geq \beta_{nh}
	\qquad \forall n \in \Omega, h \in \mathcal{H}_t
	\end{split}
	\end{equation}
	\begin{equation}
	\begin{split}
	\label{4.12}
	\alpha_{mh} - \theta_{mh}
	= \sum\limits_{n \in \Omega}\left( \gamma_{nm(h-1)} - \gamma_{nmh}
	\right) 
	\qquad \forall m \in \Omega, h \in \mathcal{H}_t
	\end{split}
	\end{equation}
	\begin{equation}
	\begin{split}
	\label{4.13}
	\sum\limits_{n \in \Omega}\left( \alpha_{nh} + \theta_{nh}
	\right) 
	\leq 1 
	\qquad \forall h \in \mathcal{H}_t
	\end{split}
	\end{equation}
	
	\paragraph{}
	The travel time constraint is formulated as equation (\ref{4.14}), where $ H_{nmh} $ is the number of time intervals required for driving and installation for the PESS to leave from one node $ n $ at time $ h $ and be prepared to operate at another node $ m $. The traveling time between the same pair of nodes may vary across time with traffic congestion. In the case studies, we estimate a time-invariant travel-time matrix for each case area by dividing the distance matrix of the area by a speed of 40 km/hour. We evaluate the impact of real traffic (modelled by time-dependent travel-time estimates) on the revenue of PESS in Figure S7. As the time-invariant travel-time matrices we applied are relatively conservative estimates in most regions and traveling time is usually short in our case study regions with a radius of 10 miles, the impact is negligible.
	\begin{equation}
	\begin{split}
	\label{4.14}
	\gamma_{nmh} \geq
	\gamma_{nm(h-1)} - \gamma_{nm(h-H_{nmh})}
	\qquad \forall n \in \Omega, m \in \Omega, h \in \mathcal{H}_t
	\end{split}
	\end{equation}
	
	\subsection{Life-cycle revenue}
	The life-cycle revenue of PESS ($ L $) in each case is calculated by aggregating all the daily revenues before the PESS life ends, as expressed in equation (\ref{6.2.1}) and (\ref{6.2.2}). The life-cycle usage or degradation limit of PESS, denoted by $ U $, is set to 2000 100\%-depth-of-discharge cycles, which is equivalently 10.8 GWh-throughput for a 2.7 MWh PESS. The price data in 2018 are repeatedly used to estimate daily revenues for each year in the PESS life.
	\begin{equation}
	\begin{split}
	\label{6.2.1}
	&L = \sum\limits_{t \in \{1,2,\dots,T\}} \delta_t \left(f(\mathbf{x}^{\ast}_{t}) +  C_{t}^{\mathrm{DEG}}(\mathbf{x}^{\ast}_{t}) \right)
	\end{split}
	\end{equation}
	\begin{equation}
	\begin{split}
	\label{6.2.2}	
	&\text{where}\ T = \max\limits_{\tau} \left\{\tau: \sum\limits_{t\in \{ 1,2,\dots,\tau\}} u_{t} \leq U\right\}	
	\end{split}
	\end{equation}
	
	\noindent where $\mathbf{x}_t^{\ast}$ is the optimizer of equation (\ref{4.1}) subject to equations (\ref{4.2})-(\ref{4.14}).
	
	\newpage
	
	\bibliographystyle{naturemag}

	
	\bibliography{Dissertation_Library}

\end{document}


\title{The economics of utility-scale portable energy storage systems in a high-renewable grid: Supplemental information}

\begingroup



\endgroup

\date{}

\maketitle

\section*{Supplemental Notes}
\subsection{Local Grid Transmission Congestion in California}
\label{SI:congestion}
\begin{figure}[!h]
	\centering
	\includegraphics[width=\textwidth]{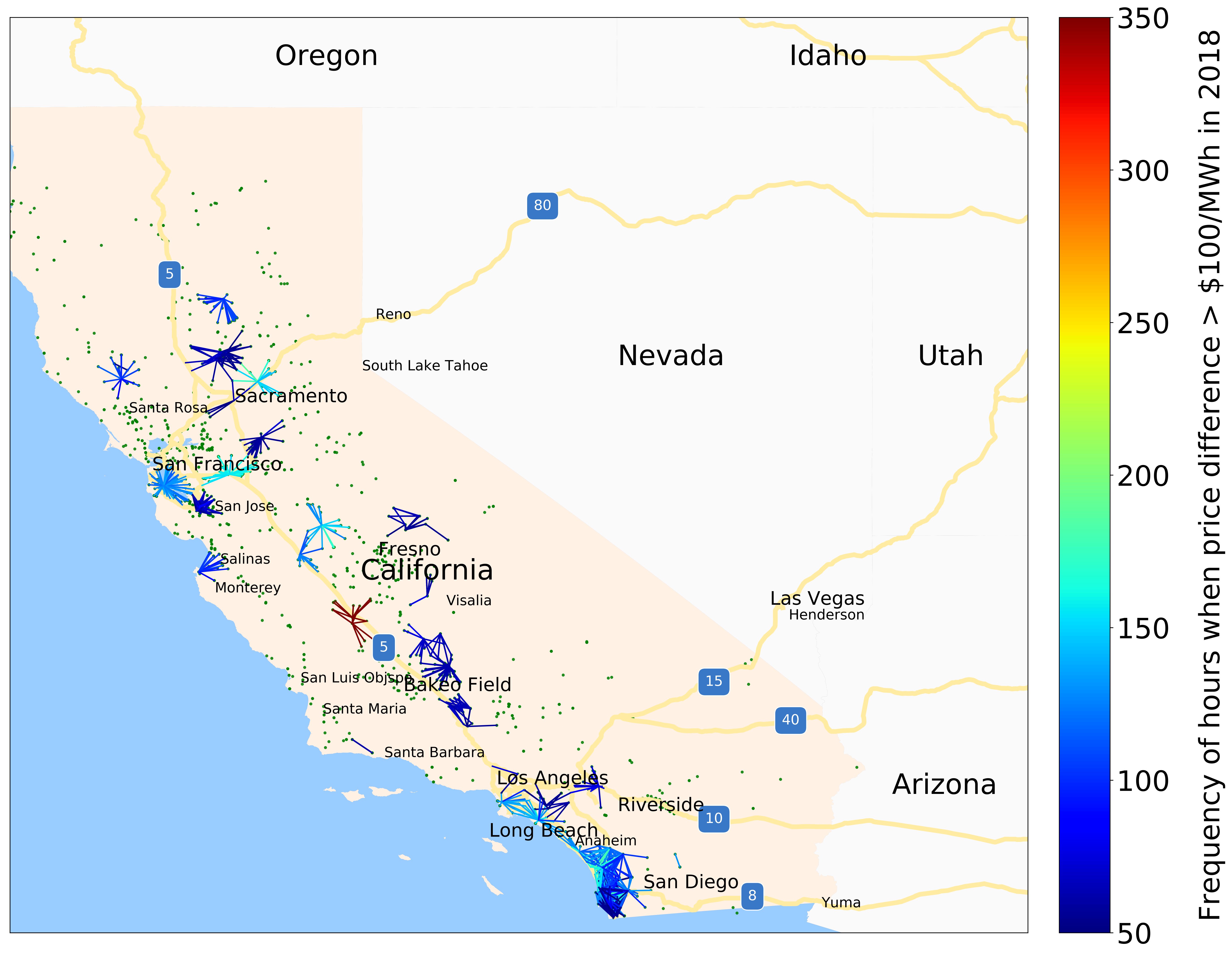}
	\caption{Distribution map of nodes with local transmission congestions in California. The green dots represent price nodes in California. The lines connect nodes that are closer than 20 miles and have more than 50 hours when the price difference is greater than \$100/MWh in 2018. The line color represents the number of hours during which the price difference is greater than \$100/MWh.}
	\label{fig-SI-1}
\end{figure}

\begin{figure}[!h]
	\centering
	\includegraphics[width=\textwidth]{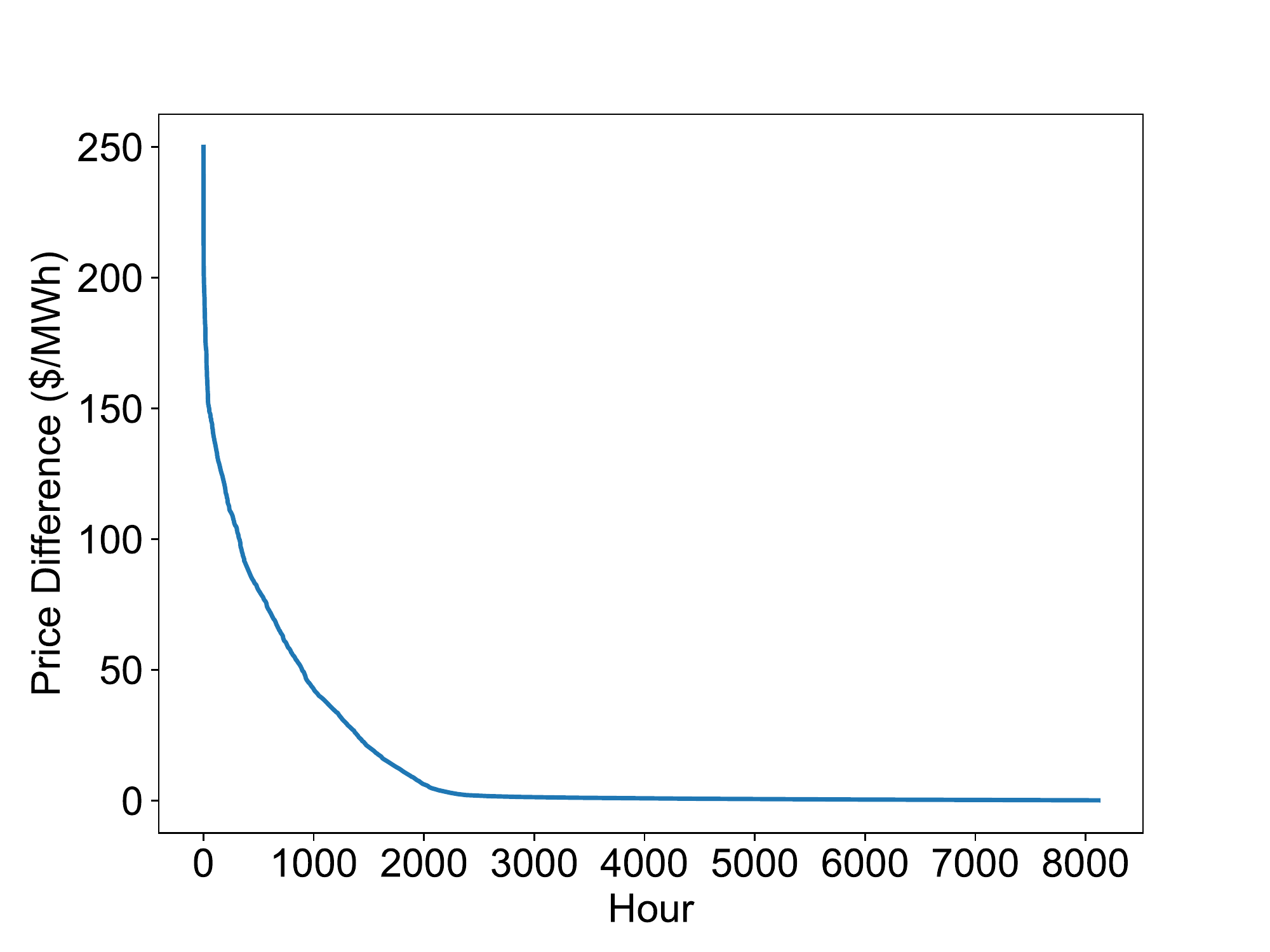}
	\caption{Frequency of price difference between two nodes around Kettleman City, CA in 2018. CAISO Node ID: KETTLEMN\_6\_N001 and HURON\_6\_N001.}
	\label{fig-SI-2}
\end{figure}

Local transmission congestions have been observed in California recently. There are significant differences in the locational marginal prices (LMPs) between some nodes with short geographical distance, and some occur very frequently. Figure \ref{fig-SI-1} presents the frequency of hours when the price difference is greater than \$100/MWh between any two nodes that are closer than 20 miles in California in 2018. \$100/MWh is a large price difference (the average LMP in California in 2018 is \$50/MWh) and comparable with the unit capital cost of lithium-ion battery packs
per cycle (\$200/kWh-capacity divided by 2000 cycles). 
\paragraph{}
  Figure \ref{fig-SI-1} shows that the network representing local congestion has a high centrality in most regions, with one or two nodes frequently having very high or low prices compared to surrounding nodes. Most congestion dyads have an annual frequency less than 200 hours. Congestions are usually more intensive around large cities.
\paragraph{}
For the two nodes with the highest frequency of large price differences (red lines in  Figure \ref{fig-SI-1}, around Kettleman City), there are over 300 hours when the arbitrage profit can be comparable with the unit-cycle capital cost of lithium-ion battery packs (\$100/MWh) and approximately 900 hours when the price difference is greater than the average LMP in California in 2018 (\$50/MWh), as shown in Figure \ref{fig-SI-2}.

\begin{figure}[!h]
	\centering
	\includegraphics[width=\textwidth]{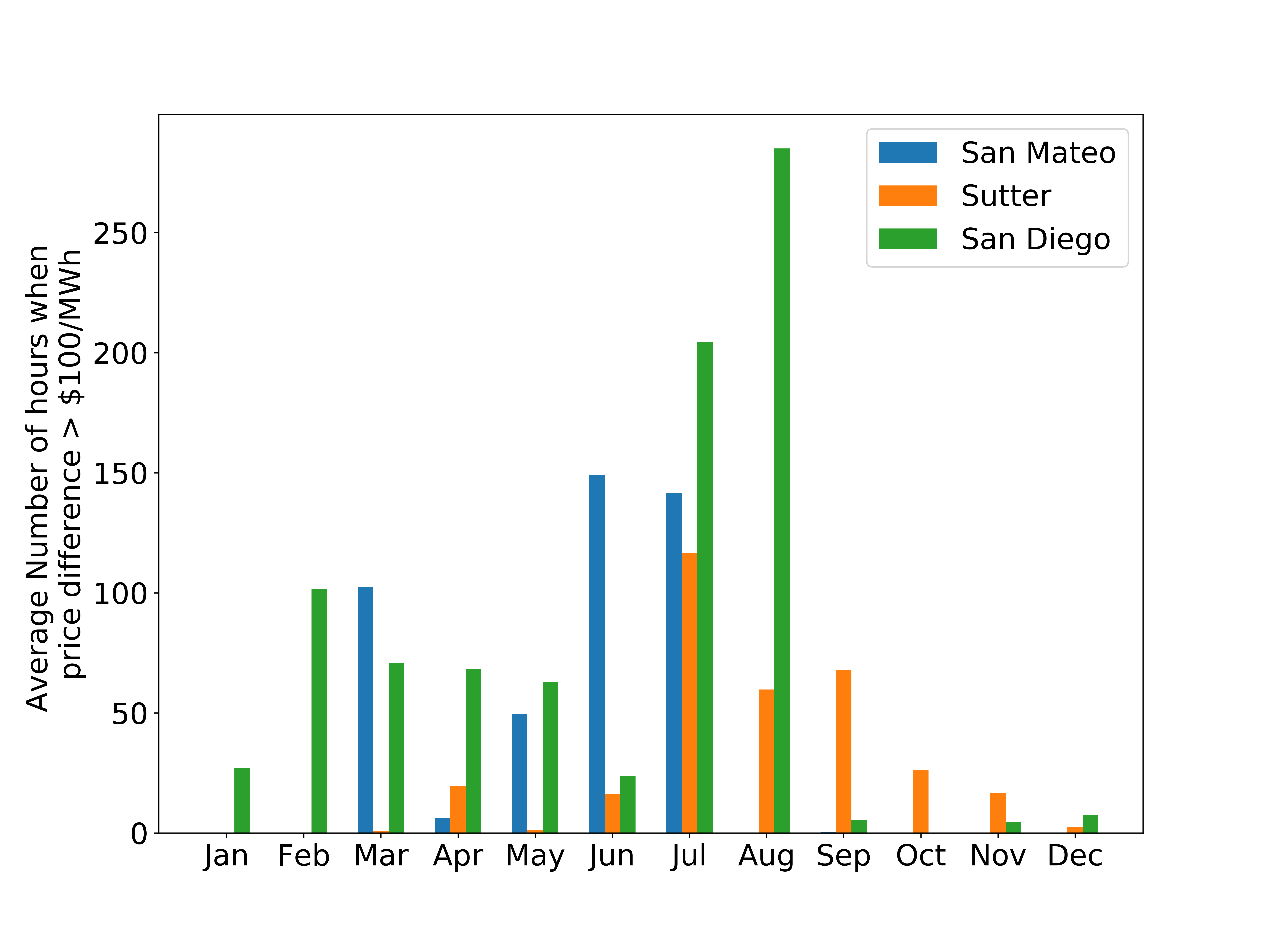}
	\caption{Monthly distributions of the average number of hours when the price difference is greater than \$100/MWh in 2018 for three counties in California, San Mateo, Sutter, and San Diego.}
	\label{fig-SI-3}
\end{figure}

\paragraph{}
Are there similar temporal patterns for congestion between different pairs of nodes? Figure \ref{fig-SI-3} presents monthly distributions of the average number of hours when the price difference is greater than \$100/MWh in 2018 for three counties in California, San Mateo, Sutter, and San Diego, and We observe no obvious temporal patterns. There are several months when more than one counties have comparatively high frequencies of significant price difference, for example, March, July, and August. However, there are also months when one county has a comparatively high frequency while the others do not, for example, February, June, and September. The correlation coefficient between the frequency series of San Mateo and San Diego in Figure \ref{fig-SI-3} is 0.19, which indicates 
a weak correlation between the congestion event frequencies. The absence of a uniform congestion pattern may favor the PESS as it can be shared among different pairs of nodes across seasons and earn higher revenues.

\paragraph{}
For transmission lines with infrequent congestion, expanding transmission capacity is also less favorable than the PESS, because the investment efficiency of the new transmission capacity will be low given its low utilization rate. In contrast, the PESS can address the congestion when needed and move to serve other locations or provide other services.

\paragraph{}
Extensive local transmission congestion could be a typical pattern for power systems in which distributed solar generation is increasing rapidly as in CAISO. Transmission line expansion usually takes a long time and cannot catch up with the rate of solar panel installation. This can also explain why the local congestion is more frequent in summer, as shown in Figure \ref{fig-SI-3}.

\subsection{The Impact of Imperfect Price Forecasts}
\label{SI:forecast}
The revenue reported in the main manuscript assumes perfect information about the day-ahead market prices of CAISO. In practice, perfect price forecasts are impossible in any electricity markets. Besides relying on advanced price forecasting technologies, one feasible strategy for PESS operators is to participate in real-time market prices while making decisions based on known day-ahead prices. Here we simulate and compare the revenue using this strategy with that assuming perfect price forecasts for a PESS travelling between the two nodes around Kettleman City with the highest frequency of large price differences in California (the same as Figure \ref{fig-SI-2}). The decision models are the same as that proposed in the Methods section in the main manuscript, and after solving the model based on day-ahead market prices, the revenue of PESS in real-time markets $ R^{\mathrm{real}}(\mathbf{x}_t) $ is calculated based on the real-time market prices $ \lambda^{\mathrm{real}}_{nh} $ and the optimal operational schedules of PESS $ P^{\mathrm{DIS}\ast}_{nh}$ and $P^{\mathrm{CHA}\ast}_{nh} $, as equation (\ref{4.s1}) below:

\begin{equation}
	\begin{split}
	\label{4.s1}
	R^{\mathrm{real}}(\mathbf{x}_t) = \Delta h\sum\limits_{h \in \mathcal{H}_t} \sum\limits_{n \in \Omega}
	\lambda^{\mathrm{real}}_{nh} \left( P^{\mathrm{DIS}\ast}_{nh} - P^{\mathrm{CHA}\ast}_{nh}
	\right)  
	\end{split}
\end{equation}

\begin{figure}[!h]
	\centering
	\includegraphics[width=\textwidth]{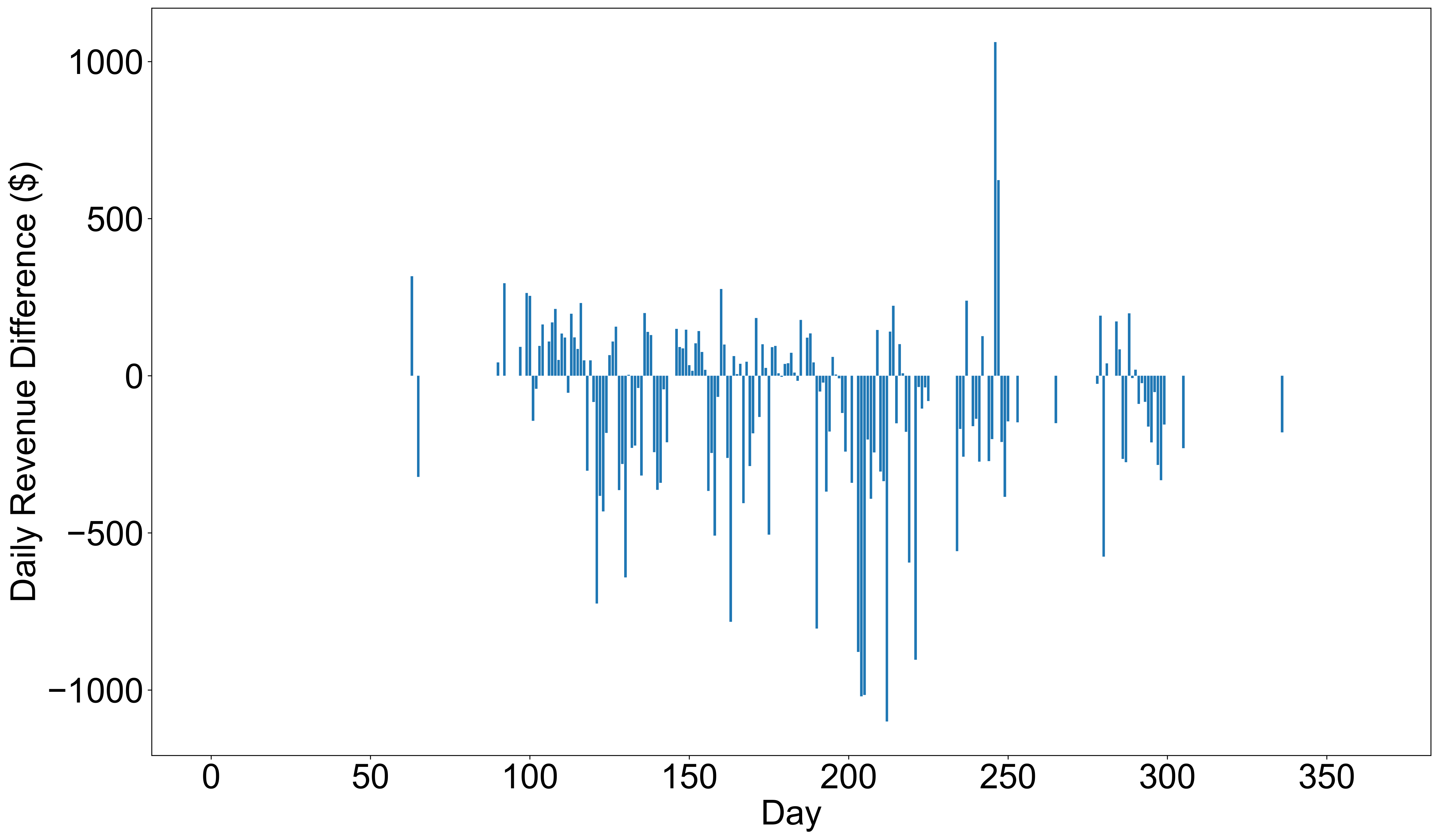}
	\caption{The differences in daily spatiotemporal arbitrage revenues between day-ahead and real-time markets in the first year of operation for the two nodes around Kettleman City. Operational and transportation schedules are optimized based on day-ahead market prices.}
	\label{fig-SI-4}
\end{figure}

\paragraph{}
The historical prices in the fifteen-minute market (FMM) in CAISO in 2018 are used in this simulation. The FMM runs approximately 37.5 minutes ahead of the dispatch interval for a horizon of 1 to 4.5 hours \cite{caiso19}. The differences in daily spatiotemporal arbitrage revenues between day-ahead and real-time markets in the first year of operation are shown in Figure \ref{fig-SI-4}. The revenue differences look random with a mean of -\$43 (real-time market revenues minus day-ahead market revenues). The total life-cycle spatiotemporal arbitrage revenue from real-time markets is approximately \$0.69 million, a 10\% revenue loss compared to the day-ahead revenues. This result is not surprising as market agents do react to the day-ahead market prices, and the transmission congestion is relieved to some degree through unit commitment adjustment in the real-time dispatch. 

\paragraph{}
Bidding in real-time markets based on day-ahead market prices is not the best strategy for the PESS, just the simplest one. The PESS operator can first bid in day-ahead markets based on price forecasts, and then adjust its bids in real-time markets to increase revenue.

\subsection{The Impact of Time-Variant Real Traffic}
\label{SI:traffic}
The travel time matrix is assumed to be constant over a year for the case studies presented in the main manuscript. As the PESS travels mostly in the daytime (Figure 5 in the main manuscript), traffic congestion is a potential problem that may delay the arrival of PESS and cause revenue loss, especially in metropolitan regions. 

\begin{figure}[!h]
	\centering
	\includegraphics[width=\textwidth]{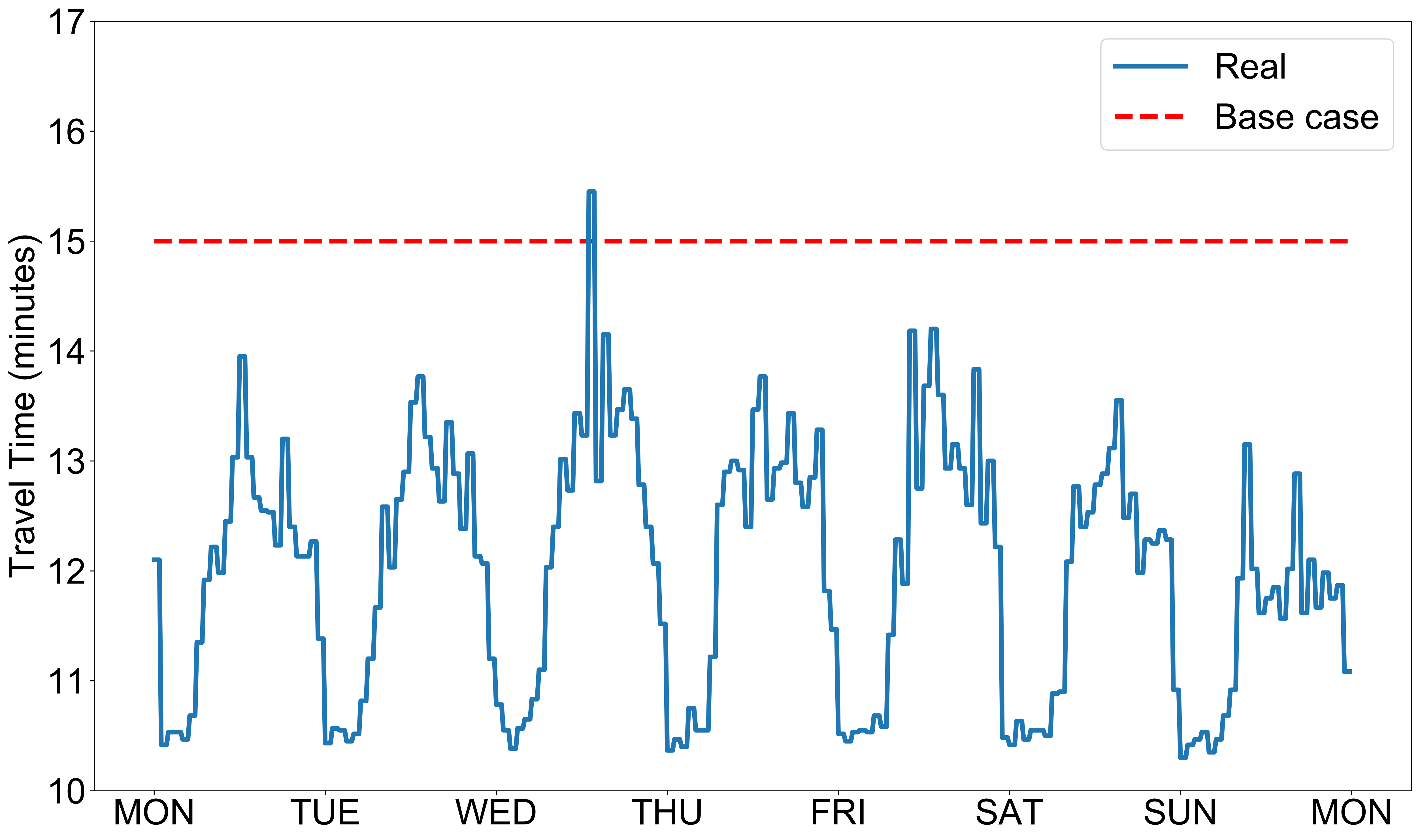}
	\caption{Weekly travel time profile for the most travelled route 24-26.}
	\label{fig-SI-5}
\end{figure}

\paragraph{}
Here we use the Microsoft Bing Maps API to estimate the realistic time-dependent travel-time matrix for the case area located in San Diego County with the highest life-cycle revenue in California. The travel-time profiles during the first week in 2018 for the most travelled route 24-26 and the top-bottom route 1-31 are presented in Figure \ref{fig-SI-5} and \ref{fig-SI-6}. We can see that the travel time of the shorter route 24-26 is less volatile than that of route 1-31. In Figure \ref{fig-SI-5}, the real travel time of route 24-26 (blue line) is mostly within 15 minutes, the constant estimate in our case studies (red dash line), and thus the traffic has little impact on the revenue estimate for this route. The same can be observed in Figure \ref{fig-SI-6}, which indicates that our constant travel-time estimate is conservative. Figure \ref{fig-SI-7} presents the daily revenue difference between our base-case model with constant travel-time matrix and the model with realistic time-dependent travel-time matrix. Both positive and negative differences are observed. The life-cycle PESS revenue difference between the two models is less than 1\%. 

\newpage

\section*{Supplemental Figures}
\label{SI:figure}
\begin{figure}[!h]
	\centering
	\includegraphics[width=\textwidth]{images/weekly_travel_time_2.png}
	\caption{Weekly travel time profile for the top-bottom route 1-31.}
	\label{fig-SI-6}
\end{figure}
\begin{figure}[!h]
	\centering
	\includegraphics[width=\textwidth]{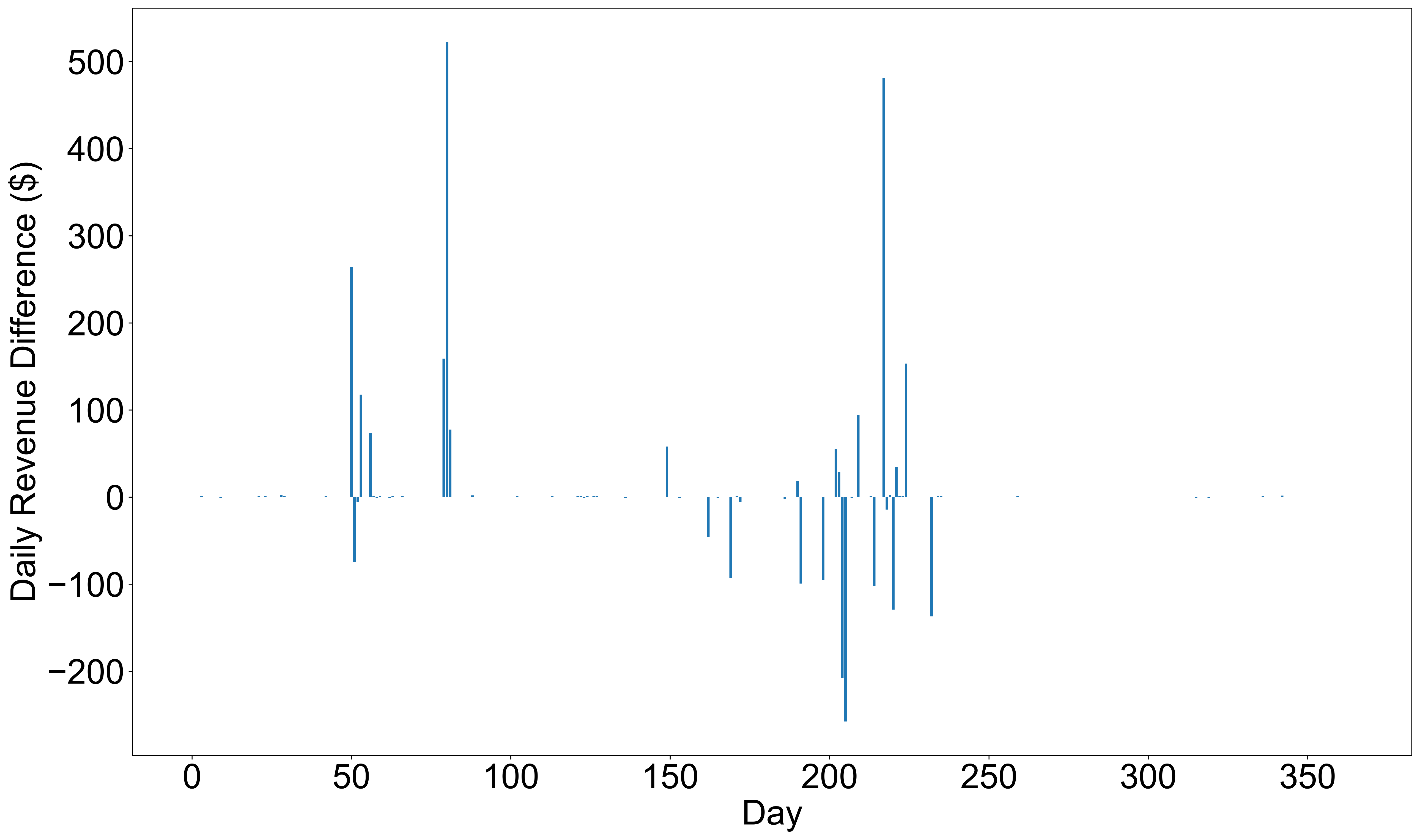}
	\caption{Daily revenue difference between models with constant and time-dependent travel-time estimates.}
	\label{fig-SI-7}
\end{figure}
\begin{figure}[!h]
	\centering
	\includegraphics[width=\textwidth]{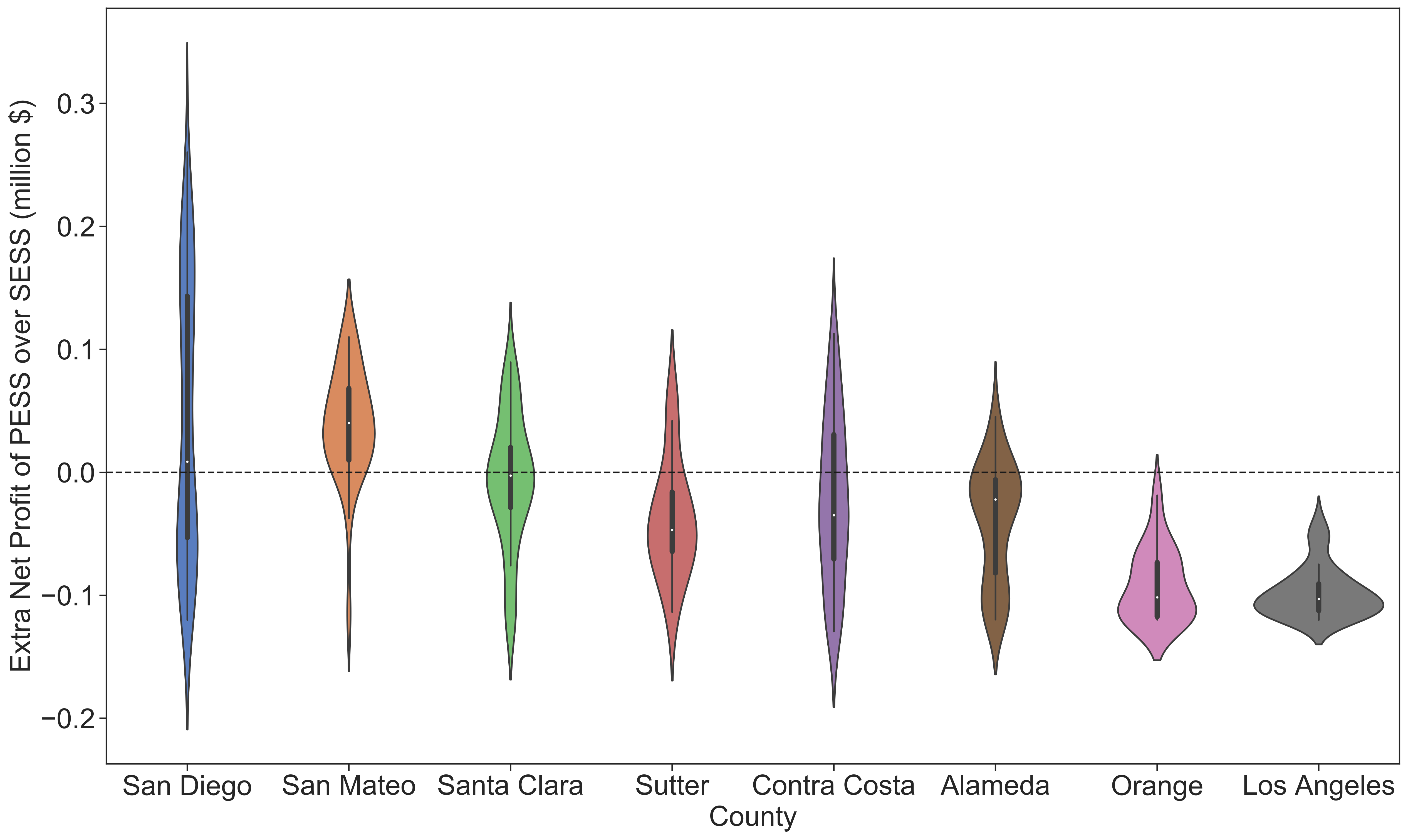}
	\caption{Distributions of the difference in estimated profitability between PESS and SESS in California by county. The top eight counties with the highest median PESS revenue and more than ten case areas are displayed and sorted in descending order of median PESS revenue. Each violin plot shows the median (the white dot), the first and third quartiles (the black bar), the upper and lower adjacent values (the black line), and the distribution of life-cycle revenues for all the case areas in each county. The up-front capital cost plus the fixed operation and maintenance (O\&M) costs are estimated based on a \$150/kWh unit cost for battery packs, while the results are similar for the case with \$200/kWh battery packs. The PESS is more profitable than the SESS in 38\% of sites in this figure and in sites covering 36\% of all 33 studied counties with at least ten case areas in California.}
	\label{fig-SI-8}
\end{figure}

\begin{figure}[!h]
		\centering
		\includegraphics[width=\textwidth]{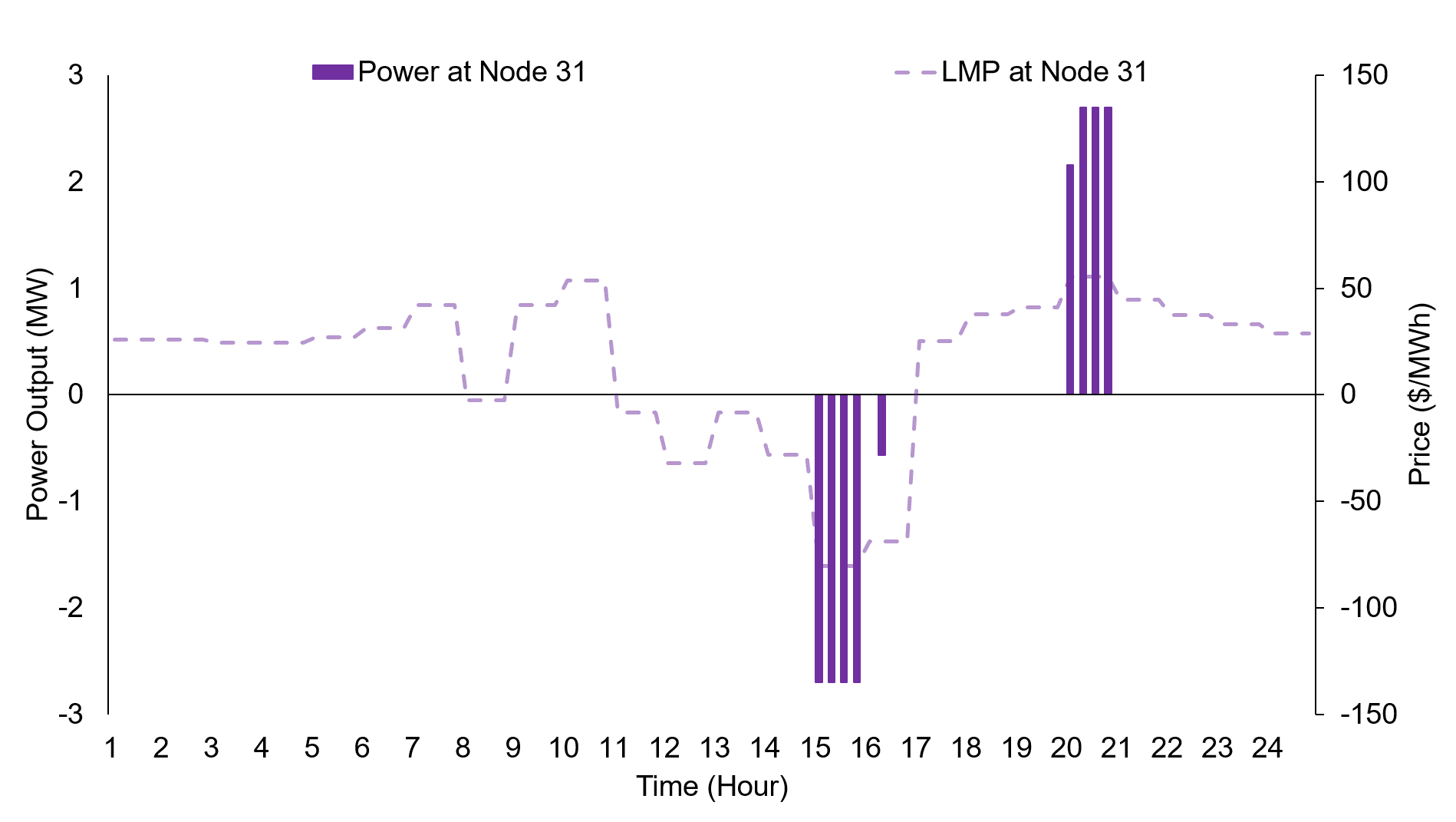}
		\caption{Optimal operational schedules for SESS doing arbitrage at node 31 in a case area in San Diego County and the locational marginal price (LMP) profiles of the node in a sample day (Mar 23, 2018). The bars indicate the charging (below the x-axis) and discharging (above the x-axis) power, and the dashed lines represent the LMP profiles.}
		\label{fig-SI-9}
\end{figure}

\clearpage


\bibliographystyle{naturemag}


\bibliography{Dissertation_Library}